# Multiplicity of a zero of an analytic function on a trajectory of a vector field

## Andrei Gabrielov

**Department of Mathematics, Purdue University**

**W. Lafayette, IN 47907-1395**

e-mail: agabriel@math.purdue.edu

Let $P(x)$ be a germ at the origin of an analytic function in $\mathbf{C}^n$, where $x = (x_1, \ldots, x_n)$, and let $\xi = \xi_1(x)\partial/\partial x_1 + \ldots + \xi_n(x)\partial/\partial x_n$ be a germ at the origin of an analytic vector field. Suppose that $\xi(0) \neq 0$, and let $\gamma$ be a trajectory of $\xi$ through the origin. Suppose that $P|_\gamma \not\equiv 0$, and let $\mu(P|_\gamma)$ be the multiplicity of a zero of $P|_\gamma$ at the origin. Let $\xi P = \xi_1 \partial P/\partial x_1 + \ldots + \xi_n \partial P/\partial x_n$ be derivative of $P$ in the direction of $\xi$, and let $\xi^k P$ be the $k$th iteration of this derivative.

We give a formula (Theorem 1) for $\mu(P|_\gamma)$ in terms of the Euler characteristic of the Milnor fibers defined by a deformation of $P, \xi P, \ldots, \xi^{n-1}P$. For a polynomial $P$ of degree $p$ and a vector field $\xi$ with polynomial coefficients of degree $q$, this allows one to compute $\mu(P|_\gamma)$ in purely algebraic terms (Theorem 2), and to give an estimate (Theorem 3) for $\mu(P|_\gamma)$ in terms of $n$, $p$, $q$, single exponential in $n$ and polynomial in $p$ and $q$. This estimate improves previous results [7, 1] which were double exponential in $n$.

For a system $\Xi = \{\xi_i\}$ of vector fields in $\mathbf{R}^n$ with polynomial coefficients of degree not exceeding $q$, this implies a single exponential in $n$ and polynomial in $q$ estimate for the degree of nonholonomy of $\Xi$, i.e., the minimal order of brackets of $\xi_i$ necessary to generate the maximal possible subspace at each point of $\mathbf{R}^n$.

For $n = 2$, our estimate coincides with the estimate for the multiplicity of a Pfaffian intersection [4, 2]. In case $n = 3$, a similar estimate was obtained in [3].

The main result of this paper can be reformulated as follows. Let $x(t) : \mathbf{C}_{t,0} \to \mathbf{C}^n_{x,0}$ be a germ of an analytic vector-function satisfying a system of nonlinear algebraic differential equations $S_i(x(t), t)dx_i/dt = Q_i(x(t), t)$ where $S_i$ and $Q_i$ are polynomial in $(x, t)$ of degree $q$, and $S_i(0, 0) \neq 0$. Let $p(t) = P(x(t), t)$ where $P$ is a polynomial in $(x, t)$ of degree $p$. Suppose that $p(t) \not\equiv 0$. Then the multiplicity of a zero of $p(t)$ at $t = 0$ can be computed in purely algebraic terms, and there is an estimate for this multiplicity in terms of $n$, $p$, $q$, single exponential in $n$ and polynomial in $p$ and $q$.

**Definition.**     Let $\epsilon \in \mathbf{C}$. A germ $\tilde{P}(x,\epsilon)$ of an analytic function at the origin in $\mathbf{C}^{n+1}$ is called a *deformation* of $P$ if $\tilde{P}(x,0) = P(x)$. For a fixed $\epsilon$, we write $\tilde{P}^\epsilon(x)$ for the function $\tilde{P}(x,\epsilon)$ considered as a function of $x$.

**Proposition 1.**     [7]. Let $\tilde{\mathbf{P}}(x,\epsilon) = \big(\tilde{P}_1(x,\epsilon),\ldots,\tilde{P}_k(x,\epsilon)\big)$ be a deformation of $\mathbf{P}(x) = \big(P_1(x),\ldots,P_k(x)\big)$. For a positive number $\delta$, let $B_\delta$ be an open ball in $\mathbf{C}^n$ of radius $\delta$ centered at the origin. Then, for a small positive $\delta$ and for a nonzero $\epsilon \in \mathbf{C}$ much smaller than $\delta$, the homotopy type of the set $\{\tilde{\mathbf{P}}^\epsilon = 0\} \cap B_\delta$ does not depend on $\delta$ and $\epsilon$. This set is called the *Milnor fiber* of the deformation $\tilde{\mathbf{P}}$.

**Proposition 2.**     [7]. Let $\tilde{\mathbf{P}}(x,\epsilon) = \big(\tilde{P}_1(x,\epsilon),\ldots,\tilde{P}_k(x,\epsilon)\big)$ be a deformation of $\mathbf{P}(x) = \big(P_1(x),\ldots,P_k(x)\big)$. Suppose that, for small $\epsilon \neq 0$, the Milnor fiber $X^\epsilon$ of $\tilde{\mathbf{P}}$ is nonsingular, i.e., $dP_1^\epsilon \wedge \ldots \wedge dP_k^\epsilon \neq 0$ at each point of $X^\epsilon$.

Let $X$ be the closure of $\bigcup_{\epsilon \neq 0} X^\epsilon$, and let $Z = X \cap \{\epsilon = 0\}$. Let $\{Z_\alpha\}$ be a good stratification (Thom's $A_\epsilon$ stratification) of $Z \setminus 0$, i.e., a Whitney stratification such that, for any sequence $(x_\nu, \epsilon_\nu)$ converging to $(x^0, 0) \in Z_\alpha$, the limit of the tangent spaces to $X^\epsilon$ at $(x_\nu, \epsilon_\nu)$, if exists, contains the tangent space to $Z_\alpha$ at $(x^0, 0)$.

Let $l(x)$ be a linear function in $\mathbf{C}^n$ such that $\{l(x) = 0\}$ is transversal to all $Z_\alpha$. Then the set of critical points of $l|_{X^\epsilon}$, for small $\epsilon \neq 0$, is zero-dimensional. Let $\nu$ be the number of these points, counted with their multiplicities, converging to the origin as $\epsilon \to 0$. Then the Milnor fiber of $\mathbf{P}$ can be obtained from the Milnor fiber of $(\mathbf{P}, l)$ by attaching $\nu$ cells of dimension $n-k$.

**Theorem 1.**     Let $P(x)$ be a germ of an analytic function in $\mathbf{C}^n$, and let $\tilde{P}(x,\epsilon)$ be a deformation of $P(x)$. Suppose that, for $k = 1,\ldots,n$, the Milnor fiber $V_k$ of $\tilde{\mathbf{P}}_k = (\tilde{P}, \xi\tilde{P}, \ldots, \xi^{k-1}\tilde{P})$ is nonsingular, and let $\chi(\tilde{\mathbf{P}}_k)$ be the Euler characteristic of $V_k$. Then

$$\mu(P|_\gamma) = \chi(\tilde{\mathbf{P}}_1) + \ldots + \chi(\tilde{\mathbf{P}}_n). \tag{1}$$

**Remark.**     A. Khovanskii suggested an alternative proof of this theorem, valid also when the Milnor fibers $V_k$ are singular, as long as $V_{n+1}$ is empty. In fact, he proved the following:

**Theorem 1′.**     Let $P(x)$ be a germ of an analytic function in $\mathbf{C}^n$, and let $\tilde{P}(x,\epsilon)$ be any deformation of $P(x)$. Let $V_k$ be the Milnor fiber of $\tilde{\mathbf{P}}_k = (\tilde{P}, \xi\tilde{P}, \ldots, \xi^{k-1}\tilde{P})$, and let $\chi(\tilde{\mathbf{P}}_k)$ be the Euler characteristic of $V_k$. Then

$$\mu(P|_\gamma) = \chi(\tilde{\mathbf{P}}_1) + \ldots + \chi(\tilde{\mathbf{P}}_\mu). \tag{2}$$

**Proof.**     Let $(z, y_2, \ldots, y_n)$ be a system of coordinates in $\mathbf{C}^n$ where $\xi = \partial/\partial z$, and let $\pi$ be projection $\mathbf{C}^n \to \mathbf{C}_y^{n-1}$. For each $y$, define $\zeta_k(y) = \chi(\pi^{-1}y \cap V_k$. Because each set $\pi^{-1}y \cap V_k$ is finite, its Euler characteristic equals to the number of points in it (not counting multiplicities). Then

$$\sum_{k=1}^{\mu} \zeta_k(y) \equiv \mu.$$

Standard "integration over Euler characteristic" arguments [9] show that

$$\int \zeta_k(y) d\chi = \chi(\tilde{\mathbf{P}}_k), \text{ and } \int \sum_{k=1}^{\mu} \zeta_k(y) d\chi = \int \mu d\chi = \mu.$$

**Proof of Theorem 1.** Let us choose a coordinate system $(z, y_2, \ldots, y_n)$ in a neighborhood of the origin in $\mathbf{C}^n$ so that $\xi = \partial/\partial z$ in this coordinate system. In particular, $\gamma$ becomes $z$-axis, and $\mu(P|_\gamma)$ equals to the multiplicity of a zero of $P(z, 0, \ldots, 0)$ at the origin.

We proceed by induction on $n$. For $n = 1$ the statement is obvious. Suppose that it holds for $n-1$, so we can apply it to the subspace $\{y_n = 0\}$ of $\mathbf{C}^n$. For a generic coordinate $y_n$, the Milnor fibers of $\tilde{\mathbf{P}}'_k = \tilde{\mathbf{P}}_k|_{y_n=0}$ are nonsingular, and the condition of Theorem 1 is satisfied. We have then

$$\mu(P|_\gamma) = \chi(\tilde{\mathbf{P}}'_1) + \ldots + \chi(\tilde{\mathbf{P}}'_{n-1}). \tag{3}$$

Let $1 \leq k \leq n - 1$. Let us fix a small nonzero $\epsilon$. According to Proposition 2, $V_k$ can be obtained from $V'_k$ by attaching $\nu_k$ cells of dimension $n - k$, where $\nu_k$ is the number of points of $V_k$ where $dy_n|_{V_k} = 0$, counted with the proper multiplicities. In particular,

$$\chi(\tilde{\mathbf{P}}'_k) = \chi(\tilde{\mathbf{P}}_k) - (-1)^{n-k} \nu_k. \tag{4}$$

The necessary transversality conditions in Proposition 2 are satisfied for a generic coordinate $y_n$, because $P(z, 0, \ldots, 0) \not\equiv 0$.

The critical points of $y_n|_{V_k}$ are defined by the linear dependence at the points of $V_k$ of the following 1-forms:

$$d(\tilde{P}^\epsilon), \ d(\xi \tilde{P}^\epsilon), \ \ldots \ d(\xi^{k-1} \tilde{P}^\epsilon), \ dy_n.$$

In other words, the rank of the following $k \times (n-1)$-matrix $A_k$ is less than $k$:

$$A_k = \begin{pmatrix} \frac{\partial}{\partial z}\tilde{P}^\epsilon & \frac{\partial}{\partial y_2}\tilde{P}^\epsilon & \ldots & \frac{\partial}{\partial y_{n-1}}\tilde{P}^\epsilon \\ \frac{\partial}{\partial z}\xi\tilde{P}^\epsilon & \frac{\partial}{\partial y_2}\xi\tilde{P}^\epsilon & \ldots & \frac{\partial}{\partial y_{n-1}}\xi\tilde{P}^\epsilon \\ \ldots & \ldots & \ldots & \ldots \\ \frac{\partial}{\partial z}\xi^{k-1}\tilde{P}^\epsilon & \frac{\partial}{\partial y_2}\xi^{k-1}\tilde{P}^\epsilon & \ldots & \frac{\partial}{\partial y_{n-1}}\xi^{k-1}\tilde{P}^\epsilon \end{pmatrix}.$$

Taking into account that $\xi = \partial/\partial z$, we find that, at the points of $X_k$, all the entries in the first column of the matrix $A_k$ are zero, except for the last entry which is $\xi^k \tilde{P}^\epsilon$.

Let $B_k$ be the matrix $A_k$ with the first column removed, and let $C_k$ be the matrix $B_k$ with the last row removed. For $k = 1, \ldots, n-2$, we have

$$B_k = C_{k+1} = \begin{pmatrix} \frac{\partial}{\partial y_2} \tilde{P}^\epsilon & \cdots & \frac{\partial}{\partial y_{n-1}} \tilde{P}^\epsilon \\ \frac{\partial}{\partial y_2} \xi \tilde{P}^\epsilon & \cdots & \frac{\partial}{\partial y_{n-1}} \xi \tilde{P}^\epsilon \\ \cdots & \cdots & \cdots \\ \frac{\partial}{\partial y_2} \xi^{k-1} \tilde{P}^\epsilon & \cdots & \frac{\partial}{\partial y_{n-1}} \xi^{k-1} \tilde{P}^\epsilon \end{pmatrix}.$$

For $k = 1, \ldots, n-1$, the rank of $A_k$ is less than $k$ if and only if either $\xi^k \tilde{P}^\epsilon = 0$ and the rank of $B_k$ is less than $k$ or $\xi^k \tilde{P}^\epsilon \neq 0$ and the rank of $C_k$ is less than $k-1$. Modifying, if necessary, $\tilde{P}$ (adding $z^k$ with a small coefficient) we can always suppose that $\xi^k \tilde{P}^\epsilon \neq 0$ at the points of $V_k$ where the rank of $C_k$ is less than $k-1$.

This means that the set of the critical points of $y_n|_{V_k}$ is a union of two disjoint sets, $X_k \cap \{\xi^k \tilde{P}^\epsilon = 0\} \cap \{\operatorname{rank} B_k < k\}$ and $X_k \cap \{\operatorname{rank} C_k < k-1\}$. Hence $\nu_k = \nu'_k + \nu''_k$, where $\nu'_k$ and $\nu''_k$ are the numbers of critical points of $y_n|_{V_k}$ in these two sets, counted with the proper multiplicities.

Taking into account that $B_k = C_{k+1}$ and $X_k \cap \{\xi^k \tilde{P}^\epsilon = 0\} = X_{k+1}$, we have $\nu'_k = \nu''_{k+1}$, for $k = 1, \ldots, n-2$. For $k = 1$, we have $\nu_1 = \nu'_1$. For $k = n-1$, we have $\nu'_{n-1} = \chi(\tilde{\mathbf{P}}_n)$, the number of points in the set $V_n$.

Replacing $\nu_k$ in (4) by $\nu'_k + \nu''_k$ and substituting (4) into (3), we see that all the values $\nu'_k$ and $\nu''_k$ cancel out, except $\nu'_{n-1}$, and (3) implies (1).

**Theorem 2.**   Let $\tilde{P}(x, \epsilon)$ be a deformation of $P(x)$ such that, for $k = 1, \ldots, n$, the Milnor fiber $V_k$ of $(\tilde{P}, \xi\tilde{P}, \ldots, \xi^{k-1}\tilde{P})$ is nonsingular. Let $l_1(x), \ldots, l_{n-1}(x)$ be generic linear forms in $\mathbf{C}^n$. For $k = 1, \ldots, n$ and $i = 1, \ldots, n-k$, let

$$\omega_{i,k} = d\tilde{P}^\epsilon \wedge d(\xi\tilde{P}^\epsilon) \wedge \ldots \wedge d(\xi^{k-1}\tilde{P}^\epsilon) \wedge dl_1 \wedge \ldots \wedge dl_{n-k-i+1} = 0.$$

Let $\nu_{0,k}$ be the number of isolated zeroes of the system

$$\tilde{P}^\epsilon = \xi\tilde{P}^\epsilon = \ldots = \xi^{k-1}\tilde{P}^\epsilon = l_1 = \ldots = l_{n-k} = 0$$

converging to the origin as $\epsilon \to 0$, and let $\nu_{i,k}$, for $i = 1, \ldots, n-k$, be the number of isolated zeroes of the system

$$\tilde{P}^\epsilon = \xi\tilde{P}^\epsilon = \ldots = \xi^{k-1}\tilde{P}^\epsilon = l_1 = \ldots = l_{n-k-i} = 0, \quad \omega_{i,k} = 0 \qquad (5)$$

converging to the origin as $\epsilon \to 0$. Then

$$\mu(P|_\gamma) = \sum_{k=1}^{n} \sum_{i=0}^{n-k} (-1)^i \nu_{i,k}. \tag{6}$$

*For a polynomial $\tilde{P}$ and a vector field $\xi$ with polynomial coefficients, equations (5) are algebraic.*

**Proof.** This theorem follows from Theorem 1 and the standard relations for the Euler characteristics of hyperplane sections (Proposition 2):

$$\chi(\tilde{P}, \xi\tilde{P}, \ldots, \xi^{k-1}\tilde{P}, l_1, \ldots, l_{n-k}) = \nu_{0,k},$$

$$\chi(\tilde{P}, \xi\tilde{P}, \ldots, \xi^{k-1}\tilde{P}, l_1, \ldots, l_{n-k-i+1}) - \chi(\tilde{P}, \xi\tilde{P}, \ldots, \xi^{k-1}\tilde{P}, l_1, \ldots, l_{n-k-i}) = (-1)^i \nu_{i,k}.$$

**Lemma 1.** *Let $l(x)$ be a germ of an analytic function such that $\xi l(0) \neq 0$. Let us choose $\delta > 0$ so that there exists a representative of $P$ in $B_\delta$. For $c = (c_1, \ldots, c_n) \in \mathbf{C}^n$, let $P_c(x) = P(x) + c_1 + c_2 l(x) + \ldots + c_n l^{n-1}(x)$.*

*(i) For a generic $c$, the set $X_{k,c} = \{P_c = \xi P_c = \ldots = \xi^{k-1} P_c = 0\} \cap B_\delta$ is nonsingular, for $k = 1, \ldots, n$.*

*(ii) For a generic $c$, the deformation $\tilde{P}(x, \epsilon) = P(x) + \epsilon(c_1 + c_2 l(x) + \ldots + c_n l^{n-1}(x))$ satisfies conditions of Theorem 1, i.e., the Milnor fiber of $(\tilde{P}, \xi\tilde{P}, \ldots, \xi^{k-1}\tilde{P})$ is nonsingular, for $k = 1, \ldots, n$.*

**Proof.** Consider coordinate system $(z, y_2, \ldots, y_n)$ where $\xi = \partial/\partial z$. Let us note first that the set $\{P_c = \xi P_c = \ldots = \xi^{k-1} P_c = 0\}$ coincides with the set

$$P_c = \frac{\partial}{\partial l} P_c = \ldots = \frac{\partial^{k-1}}{\partial l^{k-1}} P_c = 0,$$

where $\partial/\partial l$ is the partial derivative in a coordinate system $(l, y_2, \ldots, y_n)$. This follows from the chain rule. Geometrically, both sets represent the points where restriction of $P_c$ to a line parallel to $z$-axis has a zero of multiplicity at least $k$.

Now our statement is a special case of Thom's transversality theorem. Let $Z_k \subset \mathbf{C}^{2n}_{x,c}$ be defined as $\{x, c : x \in B_\delta, P_c(x) = \xi P_c(x) = \ldots = \xi^{k-1} P_c(x) = 0\}$, for $k = 1, \ldots, n$. Each set $Z_k$ is nonsingular. Let $\pi : \mathbf{C}^{2n}_{x,c} \to \mathbf{C}^n_c$ be a natural projection. The set $X_{k,c}$ is nonsingular if and only if $c$ is not a critical value of the restriction of $\pi$ to $Z_k$. Due to Sard's theorem, this holds for a general $c$.

To prove (ii), note that the set of those values of $c$ for which $X_{k,c}$ is nonsingular is a real subanalytic set which admits a complex-analytic stratification. ¿From (i), this set does not contain open subsets, hence its stratification does not have any strata of the complex dimension $n$. Hence the intersection of this set with a generic complex line through the origin is zero-dimensional.

**Theorem 3.** Let $P$ be a polynomial of degree not exceeding $p \geq n-1$, and let $\xi$ be a vector field with polynomial coefficients of degree not exceeding $q \geq 1$. Then $\mu(P|_\gamma)$ does not exceed

$$2^{2n-1} \sum_{k=1}^{n} [p + (k-1)(q-1)]^{2n}. \tag{7}$$

**Proof.** ¿From Lemma 1, there exists a deformation $\tilde{P}$ of $P$ satisfying conditions of Theorem 1, such that $P^\epsilon$ is a polynomial of degree not exceeding $p$. Hence degree of $\xi^k \tilde{P}^\epsilon$ does not exceed $p + k(q-1)$. ¿From [6], the sum of Betti numbers of the Milnor fiber $V_k$ of $(\tilde{P}, \xi\tilde{P}, \ldots, \xi^{k-1}\tilde{P})$ does not exceed $(p + (k-1)(q-1))[2p + 2(k-1)(q-1) - 1]^{2n-1}$, which does not exceed $2^{2n-1}[p + (k-1)(q-1)]^{2n}$. The estimate (7) follows now from Theorem 1.

**Theorem 4.** Let $P$ be a polynomial of degree not exceeding $p \geq n-1$, and let $\xi$ be a vector field with polynomial coefficients of degree not exceeding $q \geq 1$. Suppose that $P|_\gamma \equiv 0$, and let $P_\nu$ be any sequence of polynomials of degree not exceeding $p$ converging to $P$. Then the number of isolated zeros of $P_\nu|_\gamma$ converging to the origin as $\nu \to \infty$ does not exceed (7).

**Proof.** This follows from Theorem 3 and the results of [10]. An alternative argument was suggested by Khovanskii. Let $\mathcal{L}$ denote the linear space of all polynomials of degree not exceeding $p$ modulo polynomials identically vanishing on $\gamma$. Let $P_\nu$ be a sequence of polynomials $P_\nu$ converging to $P$ such that $M$ zeroes of $P_\nu|_\gamma$ converge to the origin as $\nu \to \infty$. These polynomials define a sequence of points $Q_\nu$ in $\mathcal{L}$. Note that the zeros of $P_\nu|_\gamma$ depend only on $Q_\nu$, and do not change when we multiply $Q_\nu$ by a constant. If we define any norm in $\mathcal{L}$, we obtain a sequence of points $Q_\nu/|Q_\nu|$ in $\mathcal{L}$ that has a non-zero limit point $Q_0$. Let $P_0$ be a polynomial of degree not exceeding $p$ such that its image in $\mathcal{L}$ is $Q_0$. Obviously, $P_0|_\gamma$ has a zero of the multiplicity $M$ at 0. Hence $M$ does not exceed (7).

**Theorem 5.** Let $\Xi = \{\xi_i\}$ be a system of vector fields in $\mathbf{C}^n$ or $\mathbf{R}^n$ with polynomial coefficients of degree not exceeding $q \geq 1$. Let $d = d(\Xi, 0)$ be dimension of the vector space $L$ spanned by the values at the origin of the vector fields $\xi_i$ and their Lie brackets of all orders. The degree of nonholonomy $N = N(\Xi, 0)$, i.e., the minimal number $N$ such that the values at the origin of the vector fields $\xi_i$ and their Lie brackets of the order not exceeding $N$ generate $L$, does not exceed

$$2^{d-2}\left(1 + 2^{2n(d-2)-2})q^{2n}\sum_{k=4}^{n+3} k^{2n}\right), \quad \text{for } d > 2, \tag{8}$$

$$1 + 2^{2n-1}q^{2n}\sum_{k=2}^{n+1} k^{2n}, \quad \text{for } d = 2. \tag{9}$$

**Proof.** According to [1], there exists a vector field $\xi$ with polynomial coefficients of degree not exceeding $2^{d-3}q$ (for $d = 2$, not exceeding $q$) and a polynomial $P$ of degree not exceeding $2^{d-1}q$ such that:
(1) $\xi(0) \neq 0$;
(2) $P|_\gamma \not\equiv 0$, where $\gamma$ is the trajectory of $\xi$ through 0;
(3) $N$ does not exceed $2^{d-2} + 2^{d-3}\mu$ (for $d = 2$, does not exceed $1 + \mu$) where $\mu$ is the multiplicity of a zero of $P|_\gamma$ at 0.

Applying this to the estimate (7) for $\mu$, we obtain (8) and (9).

**Definition 2.** (Khovanskii, unpublished; see [8].) A *Noetherian chain* of order $m$ and degree $\alpha$ is a system $f(x) = (f_1(x), \ldots, f_m(x))$ of germs of analytic functions at the origin $\mathbf{0}$ of a complex or real $n$-dimensional space, satisfying

$$\frac{\partial f_i}{\partial x_j} = g_{ij}(x, f(x)), \text{ for } i = 1, \ldots, m \text{ and } j = 1, \ldots, n, \tag{10}$$

where $g_{ij}$ are polynomials in $x$ and $f$ of degree not exceeding $\alpha \geq 1$. A function $\phi(x) = P(x, f(x))$, where $P$ is a polynomial in $x$ and $f$ of degree not exceeding $p$, is called a *Noetherian function* of degree $p$, with the Noetherian chain $f$.

The following two theorems can be reduced to Theorems 3 and 5 by adding $m$ new variables corresponding to the functions of the Noetherian chain (see [1]).

**Theorem 6.** Let $f = (f_1, \ldots, f_m)$ be a Noetherian chain of order $m$ and degree $\alpha$, and let $\xi = \sum_j \phi_j(x)\partial/\partial x_j$ be a vector field with the coefficients $\phi_j$ Noetherian of degree $q$, with the Noetherian chain $f$. Let $\psi$ be a Noetherian function of degree $p$, with the Noetherian chain $f$. Suppose that $\xi(0) \neq 0$ and that $\psi$ does not vanish identically on the trajectory $\gamma$ of $\xi$ through 0. Then the multiplicity of the zero of $\psi|_\gamma$ at 0 does not exceed

$$2^{2(n+m)-1}\sum_{k=1}^{n+m} [p + (k-1)(q+\alpha-1)]^{2(n+m)}. \tag{11}$$

**Theorem 7.**     Let $f = (f_1, \ldots, f_m)$ be a Noetherian chain in $\mathbf{C}^n$ or $\mathbf{R}^n$ of order $m$ and degree $\alpha \geq 1$. Let $\Xi = \{\xi_i\}$ be a set of vector fields with Noetherian coefficients:

$$\xi_i = \sum_j Q_{ij}(x, f(x)) \frac{\partial}{\partial x_j}$$

with $Q_{ij}$ polynomial in $x$ and $f$ of degree not exceeding $q \geq 1$. Let $d = d(\Xi, 0)$ be dimension of the vector space spanned by the values at the origin of the vector fields $\xi_i$ and their Lie brackets of all orders. The degree of nonholonomy $N(\Xi, 0)$ does not exceed

$$2^{d-2}\left(1 + 2^{2(n+m)(d-2)-2})(q+\alpha)^{2(n+m)} \sum_{k=4}^{n+m+3} k^{2(n+m)}\right), \quad \text{for } d > 2, \qquad (12)$$

$$1 + 2^{2(n+m)-1}(q+\alpha)^{2(n+m)} \sum_{k=2}^{n+m+1} k^{2(n+m)}, \quad \text{for } d = 2. \qquad (13)$$

**Remark.**     The "integration over Euler characteristic" arguments allow one to obtain an effective estimate on the multiplicity of an isolated intersection defined by Noetherian functions of degree $p$ in $n$ variables, with a Noetherian chain of order $m$ and degree $\alpha$, in terms of $n$, $m$, $\alpha$, and $p$. The proof will appear in a joint paper of A. Khovanskii and the author.

**Acknowledgements.**     This research was partially supported by NSF Grant # DMS-9704745. Part of this work was done during the visit of the author to the Fields Institute for Research in Mathematical Sciences, Toronto, Canada. The author thanks A. Khovanskii for fruitful discussions.